\documentclass[a4paper,UKenglish]{lipics_md}
\usepackage{microtype}
\usepackage{mathtools}
\usepackage[normalem]{ulem}

\usepackage{amsmath}
\usepackage{amsfonts}
\usepackage{arydshln,leftidx,mathtools}
\title{Integrating by Spheres: Summary of Blaschke--Petkantschin Formulas\footnote{
  This project has received funding from the European Research Council (ERC)
  under the European Union's Horizon 2020 research and innovation programme
  (grant agreement No 78818 Alpha).}}
\titlerunning{Integrating by Spheres}

\author[]{Anton Nikitenko}
\affil[]{IST Austria (Institute of Science and Technology Austria),
  Am Campus 1, \\ 3400 Klosterneuburg, Austria,
  \texttt{anton.nikitenko@ist.ac.at}}
\authorrunning{A. Nikitenko}

\mathsubjclass{28A75 Length, area, volume, other geometric measure theory, 60D05 Geometric probability and stochastic geometry.} 

\newcommand{\mm}[1] {\ifmmode{#1}\else{\mbox{\(#1\)}}\fi}

\newcommand{\ourproof}{\begin{proof}}
\newcommand{\eop}{\end{proof}}  %

\newcommand{\Rspace}        {\mm{{\mathbb R}}}
\newcommand{\Sspace}        {\mm{{\mathbb S}}}

\newcommand{\LGrass}[2]     {\mm{{\cal L}_{#1}^{#2}}}

\newcommand{\plane}         {\mm{P}}

\newcommand{\ERad}          {\mm{r}}

\newcommand{\uuu}           {\mm{{\bf u}}}

\newcommand{\bbb}           {\mm{{\bf b}}}
\newcommand{\xxx}           {\mm{{\bf x}}}
\newcommand{\zzz}           {\mm{{\bf z}}}
\newcommand{\yyy}           {\mm{{\bf y}}}
\newcommand{\vvv}           {\mm{{\bf v}}}
\newcommand{\www}           {\mm{{\bf w}}}

\newcommand{\Vol}[2]        {\mm{\rm Vol}_{#1}{({#2})}}

\newcommand{\abs}[1]        {\mm{\rm abs\,}{#1}}

\newcommand{\diff}          {\mm{\rm \,d}}
\newcommand{\norm}[1]       {\mm{\|{#1}\|}}

\newcommand{\Skip}[1]       {}

\begin{document}
\maketitle

\begin{abstract}
	In some applications, like some areas in stochastic geometry, a convenient change
	of variables involves spheres. In this review we summarize formulas of
	Blaschke--Petkantschin type, that help to pass from integration over $k$-tuples
	of points in space to the integration by some sphere passing through those points.
\end{abstract}

\section{Introduction}
\label{sec:intro}
Classical Blaschke--Petkantschin formulas are used to do the change of variables
in an integral over $k$-tuples of points in the space, for which only the relative
position of the points inside the subspace that they span, affinely or linearly, matters.
The following is, for example, the linear Blaschke--Petkantschin formula:
\begin{align}
 \int\displaylimits_{\xxx \in (\Rspace^n)^k} f(\xxx) \diff \xxx &= 
	\int\displaylimits_{L \in \LGrass{k}{n}}
	\int\displaylimits_{\uuu \in L^k}
		f(\uuu) \left[k! \Vol{k}{0, \uuu} \right]^{n-k} \diff \uuu \diff L.
 \label{eq:classic-bp}
\end{align}
We use this formula to explain our notation. With $\xxx$ we denote the tuple of points in $\Rspace^n$,
in this case $\xxx = (x_1, x_2, \ldots, x_k)$, and $f$ is a \textit{nice} (see below) measurable function $(\Rspace^n)^k \to \Rspace$.
The integration on the left happens with respect to $\lambda_n^{\otimes k}$ --- the $k$-fold of the standard Lebesgue measure on $\Rspace^n$.
To understand the right hand side, recall the definition of the Grassmanian $\LGrass{k}{n}$, 
the $k(n-k)$-dimensional manifold of linear $k$-planes in $\Rspace^n$. The Grassmanian has a standard associated measure,
see \cite{Had57, ScWe08}. In most applications, the function $f$ does not depend on the particular plane $L$, so
it usually integrates out to the total measure of the Grassmanian. We normalize it as
$\norm{\LGrass{k}{n}}
		= \tfrac{\sigma_n \sigma_{n-1} \cdot \ldots \cdot \sigma_{n-k+1}}
						{\sigma_1 \sigma_2     \cdot \ldots \cdot \sigma_k}$,
where $\sigma_n$ denotes the surface area of the unit ball in $\Rspace^n$.
After fixing a $k$-plane $L \subseteq \Rspace^n$ in the first integration on the right hand side, 
the inner integral goes over all $k$-tuples of points $\uuu = (u_1, u_2, \ldots u_k)$ inside $L$, with the $k$-fold of the standard Lebesgue measure on $L$, $\lambda_k^{\otimes k}$.
As $L$ is embedded in $\Rspace^n$, $f(\uuu)$ is defined, and in most applications the function does not depend on the embedding:
indeed, to understand the practical sense of the formula, one can think of a canonical example,
when $f$ is a function of the radius of the set $\xxx$. The last term in the formula is the Jacobian itself,
and $\Vol{k}{0, \uuu}$ stands for the $k$-volume of the $k$-simplex spanned by $0$ and $\uuu$.
The class of \textit{nice} functions, for which the change of variables applies, is in general quite rich, but it
is not the focus of this review to classify it: our main interest here is the Jacobian of the transformation,
so throughout the text we will just assume that $f$ is non-negative. Clearly, it implies that the theorems also hold for 
absolutely integrable functions.

As the spherical formulas we are going to talk about are conveniently used for the affine subspaces as well, we shall use the opportunity to
keep all things together, and that is why we also present the affine Blaschke--Petkantschin formula, in the form \cite[Equation (27)]{Mil71}.
It uses the fact that an affine subspace $A$ can uniquely be represented as $A=h+L$, where $L \in \LGrass{k}{n}$ and $h$ is orthogonal to $L$.
\begin{align}
	\int\displaylimits_{\xxx \in (\Rspace^n)^{k+1}} f(\xxx) \diff \xxx
		 &=  \int\displaylimits_{L \in \LGrass{k}{n}}
				 \int\displaylimits_{h \in L^\perp}
				 \int\displaylimits_{\uuu \in L^{k+1}} 
				 f(h+\uuu) \left[k! \Vol{k}{\uuu}\right]^{n-k} \diff \uuu \diff h \diff L.
  \label{eq:affine-bp}
\end{align}

The statement of \eqref{eq:classic-bp} is adapted from \cite[Section 7.2]{ScWe08}, which contains an excellent summary of
all kinds of the affine and linear formulas of this type,
and originally these results go back to works by Petkantschin \cite{Pet35}, Blaschke \cite{Bla35} and Varga \cite{Var36}.
Those formulas have served their duty well for some time, but during our recent work on discrete Morse theory
of Poisson--Delaunay mosaics, we had to find several analogues, where the integration on
the right should go over the spheres passing through the $k$ points, instead of the subspaces. To the best of author's
knowledge, the only formula explicitly stated in the literature is Theorem 7.3.1 in \cite{ScWe08}, which is stated as
Corollary \ref{col:bp_full} in this review, although
some other special cases were used in some works of Miles \cite{Mil70, Mil71}. 
With this review, we collect the new formulas
from \cite{EN17_1, EN17_2, ENR16} and prove one more.

\section{Formulas}
\label{sec:formulas}
\subsection{\texorpdfstring{$(k-1)$}{(k-1)}-sphere spanned by \texorpdfstring{$k+1$}{k+1} points in \texorpdfstring{$\Rspace^n$}{Euclidean space}}
\label{sec:BpSp}
In $\Rspace^n$, almost every $k + 1\leq n + 1$ points determine a unique $(k-1)$-sphere, which passes through them.
The opposite is also true: any $(k-1)$-sphere in $\Rspace^n$, defined by its center $z$, its radius $r$ and the linear $k$-subspace $L$, such that the sphere lives inside $z+L$,
with $k+1$ points $\uuu = (u_0, \ldots, u_k)$ chosen on it, uniquely determines a $(k+1)$-tuple of points $\xxx = (x_0, \ldots, x_k)$ in $\Rspace^n$ as
$\xxx = (z + r u_0, \ldots, z + r u_k)$, which we shorten as $\xxx = z + r \uuu$.
The following formula deals with this reparametrization, and generalizes Theorem 7.3.1 in \cite{ScWe08}.
It is proved in \cite{ENR16} for computing the expected number of simplices in a Poisson--Delaunay mosaic,
whose smallest circumspheres are empty.
	\begin{theorem}[Blaschke--Petkantschin formula for circumscribed spheres]
    \label{thm:BP-spherical}
      Let $1 \leq k \leq n$ and
      write $S(L)$ for the $(k-1)$-dimensional unit sphere
      inside $L \in \LGrass{k}{n}$, equipped with a standard measure.
      Then, for every non-negative function $f$ of $k+1$ points in $\Rspace^n$, we have
      \begin{align*}
        \int\displaylimits_{\xxx \in (\Rspace^n)^{k+1}} \!\!\! f(\xxx) \diff \xxx  &=
          \int\displaylimits_{z \in \Rspace^n}
          \int\displaylimits_{L \in \LGrass{k}{n}}
          \int\displaylimits_{r \geq 0}
          \int\displaylimits_{\uuu \in S(L)^{k+1}} \!\!\!
            f(z + r \uuu) r^{nk-1} \left[k! \Vol{k}{\uuu}\right]^{n-k+1}
                          \diff \uuu \diff r \diff L \diff z .
      \end{align*}
  \end{theorem}

The theorem can be proved by combining the affine Blaschke-Petkantschin formula \eqref{eq:affine-bp}
and the mentioned Theorem 7.3.1 from \cite{ScWe08}, which is a special case when $k=n$:
  \begin{corollary}[Top-dimensional spherical Blaschke--Petkantschin formula]
      \label{col:bp_full}
      Every measurable non-negative function
      $f \colon (\Rspace^n)^{n+1} \to \Rspace$ satisfies
      \begin{align*}
        \int\displaylimits_{\xxx \in (\Rspace^n)^{n+1}} \!\!\! f(\xxx) \diff \xxx
          &=  \int\displaylimits_{z \in \Rspace^n}
              \int\displaylimits_{r \geq 0}
              \int\displaylimits_{\uuu \in (\Sspace^{n-1})^{n+1}} \!\!\!
                f(z+r\uuu) r^{n^2-1} n! \Vol{n}{\uuu}
                \diff \uuu \diff r \diff z ,
      \end{align*}
      in which we use the standard spherical measure on the unit sphere $\Sspace^{n-1}$.
  \end{corollary}

\subsection{\texorpdfstring{$(m+q-1)$}{(m+q-1)}-sphere spanned by \texorpdfstring{$m$}{m} points in \texorpdfstring{$\Rspace^n$}{Euclidean space} and \texorpdfstring{$q+1$}{q+1} fixed points or a \texorpdfstring{$(q-1)$}{(q-1)}-circle}
\label{sec:fixed}
The formulas from this section appear to be new, and we give a proof for them in Section \ref{sec:fixed-proof}.
We start with a special case.
Notice that for almost every $n$ points $\xxx = (x_1, \ldots, x_n)$ in $\Rspace^n$ there
exists a unique $(n-1)$-sphere passing through $\xxx$ and $0$.
If $Z$ is the center of this sphere, the radius
is $r = \norm{Z}$, because the sphere passes through $0$, so we rather choose $r \in (0, +\infty)$
and $z$ as a point on the unit sphere $\Sspace^{n-1}$ as new parameters, to make $Z = r z$.
We thus use the following reparametrization:
take $r$ as the radius, take $z \in \Sspace^{n-1}$ as the direction of the center,
and take $\uuu \in (\Sspace^{n-1})^n$ for the points, to have
$\xxx = r z + r \uuu$. The new formula then reads:
	
	\begin{theorem}[First Blaschke--Petkantschin formula for pivoted spheres]
    \label{thm:BP-pivotal-1}
      For every non-negative function $f$ of $n$ points in $\Rspace^n$ we have
      \begin{align*}
        \int\displaylimits_{\xxx \in (\Rspace^n)^n} \!\!\! f(\xxx) \diff \xxx  &=
          \int\displaylimits_{r \geq 0}
          \int\displaylimits_{z \in \Sspace^{n-1}}
          \int\displaylimits_{\uuu \in (\Sspace^{n-1})^{n}} \!\!\!
            f(r z + r \uuu) r^{n^2-1} n! \Vol{n}{-z, \uuu}
                          \diff \uuu \diff z \diff r.
      \end{align*}
  \end{theorem}	
	
Combining this with the linear Blaschke--Petkantschin formula \eqref{eq:classic-bp}, we can
easily get an equation for $m$ points in $\Rspace^n$:
	\begin{theorem}[Second Blaschke--Petkantschin formula for pivoted spheres]
    \label{thm:BP-pivotal-2}
      For every non-negative function $f$ of $m$ points in $\Rspace^n$ we have
      \begin{align*}
        \MoveEqLeft\!\!\!\!\!\!\!\int\displaylimits_{\xxx \in (\Rspace^n)^m} \!\!\!\!\!\!\! f(\xxx) \diff \xxx  = \!\!\!
				  \int\displaylimits_{L \in \LGrass{m}{n}}
				  \int\displaylimits_{r \geq 0}
          \int\displaylimits_{z \in S(L)} 
				  \int\displaylimits_{\uuu \in (S(L))^{m}} \!\!\!\!\!\!\!\!
            f(r z \! + \! r \uuu) r^{m n-1} \left[ m! \Vol{m}{-z, \uuu} \right]^{n-m+1}
                          \diff \uuu \diff z \diff r \diff L,
      \end{align*}
			where $S(L)$ is the unit sphere in $L$.
  \end{theorem}	

For the generic statement, consider $q+1$ points $p_i$ in general position in $\Rspace^n$.
Then, for almost any $m \leq n-q$ points $x_j$, there exists a unique $(m+q-1)$-sphere passing through
all of $x_j$ and $p_i$.
Assume, without loss of generality, that the center of the $(q-1)$-sphere, circumscribed around $p_i$, is the origin.
Denote by $r_0 \geq 0$ the radius of this sphere, and by $Q$ the unique $q$-hyperplane that contains all $p_i$.
An $(m+q-1)$-sphere of radius $r$ passes through these points iff its center lies in $Q^\perp$, the $(n-q)$-subspace of $\Rspace^n$
orthogonal to $Q$, and 
the distance from the center to the origin is $\sqrt{r^2-r_0^2}$.
In other words, the question is equivalent to finding an $(m+q-1)$-sphere that contains a fixed
$(q-1)$-sphere in $Q$ of radius $r_0$ centered at the origin, $S(0, r_0, Q)$,
and passes through all $x_j$.
The reparametrization in question is thus:
take an $m$-dimensional hyperplane $L$ orthogonal to $Q$, so that all points lie in $L \oplus Q$,
take radius $r \geq r_0$,
take a point $z$ on the unit $(m-1)$-sphere $S(L)$ (thus letting the center of the circumsphere be $\sqrt{r^2-r_0^2} z$),
take $m$ points $\uuu$ on the unit $(m+q-1)$-sphere $S(L \oplus Q)$,
and let $\xxx = \sqrt{r^2-r_0^2} z + r \uuu$. Then the following transformation can be applied:
	\begin{theorem}[Blaschke--Petkantschin formula for spheres containing a circle]
    \label{thm:BP-pivotal-3}
			Fix a $q$-plane $Q \in \LGrass{q}{n}$ and radius $r_0 \geq 0$, and let $m \leq n - q$.
      Then, for every non-negative function $f$ of $m$ points in $\Rspace^n$, we have
      \begin{align*}
        \int\displaylimits_{\xxx \in (\Rspace^n)^m} \!\!\! f(\xxx) \diff \xxx  &= \!\!
				  \int\displaylimits_{L \in \LGrass{m}{n-q}}
				  \int\displaylimits_{r \geq r_0}
          \int\displaylimits_{z \in S(L)} 
				  \int\displaylimits_{\uuu \in (S(L \oplus Q))^{m}} \!\!\!      
						f(\sqrt{r^2-r_0^2} z + r \uuu) \\
						 &~~~r^{m(n-1)+1} (r^2-r_0^2)^{\tfrac{m-2}{2}}
						 \left[ m! \Vol{m}{-\tfrac{\sqrt{r^2-r_0^2}}{r}\zzz, \uuu^{L}}\right]^{n-q-m+1}
														\!\!\!\diff \uuu \diff z \diff r \diff L,
      \end{align*}
			where \LGrass{m}{n-q} is the Grassmanian of $m$-hyperplanes orthogonal to $Q$.
			$S(L)$ is the unit sphere in $L$, whereas $S(L \oplus Q)$ is the unit sphere in $L \oplus Q$,
			and $\uuu^{L}$ is the orthogonal projection of $\uuu$ onto $L$.			
  \end{theorem}	
\subsection{\texorpdfstring{$k+1$}{k+1} points in \texorpdfstring{$\Rspace^n$}{Euclidean space} anchored at \texorpdfstring{$\Rspace^m$}{a subspace}}
\label{sec:Anchored}
    The next formula was proved during the analysis of weighted Poisson-Delaunay mosaics in \cite{EN17_2}.
		Fix an $m$-dimensional affine subspace $F_m$ in $\Rspace^n$.
		We are interested in integrating by spheres centered (aka \emph{anchored})
		in this subspace: for almost every $k+1 \leq m+1$ points $\xxx$ in $\Rspace^n$,
		there is a unique smallest sphere, centered at a point in $F_m$
		passing through $\xxx$. It is not hard to see, that the center of this sphere
		lies in the orthogonal projection $P_k$ onto $F_m$ of the affine $k$-plane passing through $\xxx$.
		The opposite mapping can be defined as follows:
		take $z \in F_m$ as the center of the sphere and $r$ as its radius,
		take $P_k$ as a $k$-dimensional plane inside $F_m$,
		and take $\uuu = (u_0, \ldots, u_k)$ on the unit $(k-1)$-sphere $S$ in $P_k \times F_m^\perp$,
		where $F_m^\perp$ is the $(n-m)$-dimensional linear subspace orthogonal to $F_m$.
		Then, $\xxx = z + r \uuu$ uniquely determines almost every $(k+1)$-tuple of points
		in $\Rspace^n$.		
		The corresponding statement follows. For $k=m$ it is Theorem \ref{thm:BP-spherical}.
    \begin{theorem}[Anchored Blaschke--Petkantschin formula]
      \label{thm:BP-spherical-anchored}
      Let $0 \leq k \leq m \leq n$ and $\alpha = n(k+1)-(m+1)$.
			Fix an affine $m$-subspace $F_m$ of $\Rspace^n$.
      Then every measurable non-negative function
      $f \colon (\Rspace^n)^{k+1} \to \Rspace$ satisfies
      \begin{align*}
        \!\!\!\!\int\displaylimits_{\xxx \in (\Rspace^n)^{k+1}}
              \!\!\!\!\!\!f(\xxx) \diff \xxx
          &=  \int\displaylimits_{z \in F_m}
              \int\displaylimits_{r \geq 0}
              \int\displaylimits_{P_k \in \LGrass{k}{m}}
              \int\displaylimits_{\uuu \in (S)^{k+1}}
                \!\!\!\!\!f(z+r\uuu) r^\alpha \left[k! \Vol{k}{\uuu^{P_k}}\right]^{m-k+1}
                \diff \uuu \diff P_k \diff r \diff z ,
      \end{align*}
      in which $\LGrass{k}{m}$ is the Grassmannian of (linear) $k$-planes
      in $F_m$, $\uuu^{P_k}$ is the projection of $\uuu$ to $P_k$,
      and $S$ is short for the unit sphere in $P_k \times F_m^\perp$.
    \end{theorem}

\subsection{\texorpdfstring{$k+1$}{k+1} points on \texorpdfstring{$\Sspace^n$}{sphere}}
\label{sec:sphere}
	The formula in this section lives in another space:
	it integrates over the circles on the $n$-sphere $\Sspace^n \subset \Rspace^{n+1}$.
	More precisely, if $k+1$ points $\xxx$ are placed on $\Sspace^n$, there is
	a unique $(k-1)$-sphere inside $\Sspace^n$ passing through $\xxx$, and we
	aim at integrating over such spheres.
	This formula was explicitly stated in \cite{EN17_1} and applied to study spherical Poisson--Delaunay mosaics,
	and a special case was implicitly used by Miles in \cite{Mil71} to compute the moments
	of volumes of isotropic random simplices.
	
	To express the result, we write $\plane^\perp$ for the $(n-k+1)$-plane
	orthogonal to the $k$-plane $\plane$,
	both passing through the origin in $\Rspace^{n+1}$,
	and we write $S_\plane$ for the unit $(k-1)$-sphere in $\plane$.
		\begin{theorem}[Blaschke--Petkantschin formula on the sphere]
      \label{thm:BPforSphere}
      Let $1 \leq k \leq n$,
      and let $f \colon (\Sspace^n)^{k+1} \to \Rspace$ be a non-negative measurable function.
      Then
      \begin{align*}
        \!\! \int\displaylimits_{\xxx \in (\Sspace^n)^{k+1}}
              \!\!\!\!\!\!\! f(\xxx) \diff \xxx
          &=  \!\!\!\! \int\displaylimits_{\plane \in \LGrass{k}{n+1}}
              \int\displaylimits_{p \in \plane^\perp} \!\! \ERad^{kn-2}
              \!\!\!\!\! \int\displaylimits_{\uuu \in (S_\plane)^{k+1}}
                \!\!\!\!\!\!\! f(p + \ERad \uuu) \left[ k!\Vol{k}{\uuu} \right]^{n-k+1}
              \diff \uuu \diff p \diff \plane ,
      \end{align*}
      in which $\ERad^2 = 1-\norm{p}^2$, implicitly assuming $\norm{p} \leq 1$,
      If $f$ is rotationally symmetric, we define
      $f_\ERad (\uuu) = f(p + \ERad \uuu)$, in which
      $\uuu$ is a $k$-simplex on $\Sspace^{k-1} \subset \Rspace^k$,
      and $p$ is any point with $\norm{p}^2 = 1 - \ERad^2 \leq 1$ in the $(n-k+1)$-plane
      orthogonal to $\Rspace^k \subseteq \Rspace^{n+1}$.
      With this notation, we have
      \begin{align*}
        \!\! &\int\displaylimits_{\xxx \in (\Sspace^n)^{k+1}}
              \!\!\!\!\! f(\xxx) \diff \xxx
          = 
            \tfrac{\sigma_{n+1}}{2} \norm{\LGrass{k}{n}}
              \int\displaylimits_{t=0}^1
                t^{\frac{kn-2}{2}}\!(1\! -\! t)^{\frac{n-k-1}{2}}
              \!\!\!\!\!\!\! \int\displaylimits_{\uuu \in (\Sspace^{k-1})^{k+1}}
                \!\!\!\!\!\!\!\! f_{\sqrt{t}} (\uuu) \left[ k!\Vol{k}{\uuu} \right]^{n-k+1}
              \! \diff \uuu \diff t .
      \end{align*}
    \end{theorem}

\section{Proof for a fixed subcircle}
\label{sec:fixed-proof}
Theorems \ref{thm:BP-pivotal-1} and \ref{thm:BP-pivotal-2} are the special cases of Theorem \ref{thm:BP-pivotal-3},
so we prove the latter. We first assume $m = n-q$, and then apply the linear formula \eqref{eq:classic-bp}
to get to arbitrary $m$. We follow the approach in \cite{Mol94} to prove the theorem. We will slightly
abuse the notation to emphasize that some objects are vectors: we write $\xxx_i=x_i$, $\zzz = z$, and $\uuu_i = u_i$.
We also write $r^* = \sqrt{r^2-r_0^2}$ and $r' = \tfrac{\partial r^*}{\partial r} = \tfrac{r}{\sqrt{r^2-r_0^2}} = \tfrac{r}{r^*}$.

If $m=n-q$, then $\LGrass{m}{n-q}$ contains only $Q^\perp$, so the change of variables turns be
$\xxx = r^* \zzz + r \uuu$ where $\zzz \in S(Q^\perp)$, and $\uuu \in (S(\Rspace^n))^{n-q}$.
We need to compute the Jacobian of this change of variables, which is the determinant
of the matrix of partial derivatives. We do it locally at each point, choosing some local coordinates for $\zzz$ and $\uuu_i$
on the corresponding spheres. Like in \cite{ScWe08}, we write the matrices of partial derivatives with respect to these
coordinates as $\dot{\zzz}$ and $\dot{\uuu}_i$. 
The dimensions of the matrices are then as following:
$\zzz$ and $\uuu_i$ are column vectors $n \times 1$,
$\dot{\zzz}$ is an $n \times (m-1)$ matrix, and $\dot{\uuu}_i$ are $n \times (n-1)$ matrices.
We can choose the local coordinates on the spheres and in $\Rspace^n$ in such a way that
the matrices $[\uuu_i \dot{\uuu}_i]$ are all orthogonal, 
$Q$ is the subspace where first $m$ coordinates are $0$, and the matrix $[\zzz \dot{\zzz}] = E_{n,m}$,
an $m \times n$ matrix with ones on the diagonals and zeros everywhere else.
With this, we can get a determinant expression for the Jacobian. Using block notation
for matrices, we get for the Jacobian $\frac{\partial(\xxx_1, \ldots, \xxx_m)}{\partial(r, \zzz, \uuu_1, \ldots, \uuu_m)}$:
        \begin{align*}
          J (r, \zzz, \uuu)
             &=  \abs{ \left|
            \begin{array}{cc:cccc}
              r' \zzz + \uuu_1 & r^* \dot{\zzz}  & r \dot{\uuu}_1 & 0                 & \ldots & 0 \\
              r' \zzz + \uuu_2 & r^* \dot{\zzz}  & 0              & r \dot{\uuu}_2    & \ldots & 0 \\
              \vdots           & \vdots          & \vdots         & \vdots            & \ddots & \vdots \\
              r' \zzz + \uuu_m & r^* \dot{\zzz}  & 0              & 0                 & \ldots & r \dot{\uuu}_m
            \end{array} \right|}.
        \end{align*}
We extract the powers of radius, and from now on we will work with $\tilde{J} = \tfrac{J}{(r^*)^{m-1} r^{m(n-1)}}$.
Since transposing the matrix does not affect the determinant, we can write:
\begin{align*}
	\MoveEqLeft\!\!\!\!\!\!\!\!\!\!\tilde{J}^2
		 =  		
		\left|
		\begin{array}{cccc}
		  r' \zzz^T + \uuu^T_1 & r' \zzz^T + \uuu^T_2 & \ldots &	r' \zzz^T + \uuu^T_m \\
		  \dot{\zzz}^T         & \dot{\zzz}^T         & \ldots & \dot{\zzz}^T          \\
			\hdashline
			\dot{\uuu}^T_1       & 0                    & \ldots & 0                     \\
			0                    & \dot{\uuu}^T_2       & \ldots & 0                     \\
			\vdots               & \vdots               & \ddots & \vdots                \\
			0                    & 0                    & \ldots & \dot{\uuu}^T_m
		\end{array} \right|
		\!%
		\left|
		\begin{array}{cc:cccc}
			r' \zzz + \uuu_1 & r^* \dot{\zzz}  & r \dot{\uuu}_1 & 0                 & \ldots & 0 \\
			r' \zzz + \uuu_2 & r^* \dot{\zzz}  & 0              & r \dot{\uuu}_2    & \ldots & 0 \\
			\vdots           & \vdots          & \vdots         & \vdots            & \ddots & \vdots \\
			r' \zzz + \uuu_m & r^* \dot{\zzz}  & 0              & 0                 & \ldots & r \dot{\uuu}_m
		\end{array} \right|\!\hspace{-0.81pt}.
\end{align*}
	The following identities are required to proceed.
	The orthogonality of the matrices $[\uuu_i \dot{\uuu}_i]$ implies that
	$\uuu_i^T \uuu_i = 1$,
	$\dot{\uuu}_i^T \dot{\uuu}_i = E_{n-1,n-1}$,
	whereas $\uuu_i^T \dot{\uuu}_i = 0_{1,n-1}$ is the zero row vector of length $n-1$.
	It further implies that $\dot{\uuu}_i \dot{\uuu}_i^T = E_{n,n} - \uuu_i \uuu_i^T$, but we will use
	this observation later.
	We also use $\zzz^T \zzz = 1$, $\dot{\zzz}^T \dot{\zzz} = E_{m-1,m-1}$, and $\zzz^T \dot{\zzz} = 0_{1,m-1}$.
	We now multiply the matrices:
\begin{align*}
	\tilde{J}^2
		 &=  		
		\left|
		\begin{array}{cc:ccc}
		  \sum{(r' \zzz^T + \uuu^T_i)(r' \zzz + \uuu_i)} & \sum{\uuu^T_i \dot{\zzz}} &	r' \zzz^T \dot{\uuu}_1   & \ldots &	r' \zzz^T \dot{\uuu}_m \\
		  \sum{\dot{\zzz}^T \uuu_i}                      & m E_{m-1,m-1}             &    \dot{z}^T \dot{\uuu}_1 & \ldots & \dot{z}^T \dot{\uuu}_m   \\
			\hdashline
			r' \dot{\uuu}^T_1 \zzz                         & \dot{\uuu}^T_1 \dot{z}    &  E_{n-1, n-1}             & \ldots & 0 \\
			\vdots                                         & \vdots                    &  \vdots                   & \ddots & \vdots \\
			r' \dot{\uuu}^T_m \zzz                         & \dot{\uuu}^T_m \dot{z}    &  0                        & \ldots & E_{n-1, n-1}
		\end{array} \right|.
\end{align*}
We want to reduce the determinant to the top left corner indicated by the dashed lines. For this we zero out the top right corner,
by subtracting the third to last block rows with appropriate coefficients. We then arrive at, before any simplifications:
\begin{align*}
	\tilde{J}^2
		 &=  		
		\left|
		\begin{array}{cc}
		  \sum{(r' \zzz^T + \uuu^T_i)(r' \zzz + \uuu_i)} - \sum{r'^2 \zzz^T \dot{\uuu}_i \dot{\uuu}^T_i \zzz}  &  \sum{\uuu^T_i \dot{\zzz}} - \sum{r' \zzz^T \dot{\uuu}_i \dot{\uuu}^T_i \dot{\zzz}} \\
			\sum{\dot{\zzz}^T \uuu_i} - \sum{r'\dot{\zzz}^T\dot{\uuu}_i\dot{\uuu}^T_i\zzz} & m E_{m-1,m-1} - \sum{\dot{\zzz}^T \dot{\uuu}_i \dot{\uuu}^T_i \dot{\zzz}}
		\end{array} \right|.
\end{align*}
We simplify this matrix, $X$, memberwise.
Recall that $\dot{\uuu}_i \dot{\uuu}_i^T = E_{n,n} - \uuu_i \uuu_i^T$.
Applying this observation, noticing that $X_{12} = X^T_{21}$, writing $\langle \cdot, \cdot \rangle$ for the inner product in $\Rspace^n$, and opening the brackets, we get:
\begin{align*}
	X_{11} &= \sum \left( r'^2 + 2 r' \langle \zzz, \uuu_i \rangle + 1 - r'^2 \zzz^T(E_{n,n} - \uuu_i \uuu^T_i) \zzz \right) \\
	       &= \sum \left(2 r' \langle \zzz, \uuu_i \rangle + 1 + r'^2 \langle \zzz, \uuu_i \rangle ^ 2\right) \\
			   &= \sum \left(\langle r' \zzz, \uuu_i \rangle + 1\right)^2; \\
	X_{12} &= X^T_{21} 
	       = \sum \left( \uuu^T_i \dot{\zzz} - r' \zzz^T (E_{n,n} - \uuu_i \uuu^T_i) \dot{\zzz} \right) \\
	       &= \sum \left( \uuu^T_i \dot{\zzz} + r' \zzz^T \uuu_i \uuu^T_i \dot{\zzz} \right) \\
				 &= \sum ( \langle r' \zzz, \uuu_i \rangle + 1) \uuu^T_i \dot{\zzz}; \\
	X_{22} &= m E_{m-1, m-1} - \sum \dot{\zzz}^T (E_{n,n} - \uuu_i \uuu^T_i) \dot{\zzz} \\
	       &= m E_{m-1, m-1} - \sum \left(E_{m-1,m-1} - \dot{\zzz}^T\uuu_i \uuu^T_i\dot{\zzz} \right) \\
				 &= \sum (\uuu^T_i\dot{\zzz})^T \uuu^T_i\dot{\zzz}.	
\end{align*}
We used that $\langle \vvv, \www \rangle = \vvv^T  \www = \www^T \vvv$. We now employ a simple lemma, verification of which is trivial:
\begin{lemma}
	Let $a_i$, $1 \leq i \leq m$, be real numbers, and 
	$\bbb_i$, $1 \leq i \leq m$, be row vectors of length $m - 1$.
	Then, using block notation, the following matrix
	can be decomposed as a product of two $m \times m$ matrices, which
	are a transposition of each other:
	\begin{align*}
		\left(
		\begin{array}{cc}
		  \sum a_i^2         &  \sum a_i \bbb_i \\
		  \sum a_i \bbb^T_i  &  \sum \bbb^T_i \bbb_i
		\end{array} 
		\right) &=
		\left(
		\begin{array}{cccc}
		  a_1      & a_2      & \ldots & a_m \\
			\bbb^T_1 & \bbb^T_2 & \ldots & \bbb^T_m
		\end{array} 
		\right)
		\cdot
		\left(
		\begin{array}{cc}
		  a_1      & \bbb_1 \\
			a_2      & \bbb_2 \\
			\vdots   & \vdots \\
			a_m      & \bbb_m.
		\end{array} 		
		\right).\tag*\qed
  \end{align*}
\end{lemma}
We apply the lemma with $a_i = \langle r' \zzz, \uuu_i \rangle + 1$ and 
$b_i = \uuu^T_i \dot{\zzz}$. Recall that $[\zzz \dot{\zzz}] = E_{n, m}$,
so $a_i = r' \uuu_{i,1} + 1$ and $b_{ij} = \uuu_{i,j+1}$, where $\uuu_{i,k}$ is the $k$-th
coordinate of $\uuu_i$. Since the transposition does not change the
determinant, we get, dividing the first row of the matrix by $r'$,
\begin{align*}
	\tilde{J}
		 &=
		r' \cdot
		\abs{
		\left|
		\begin{array}{cccc}
		  u_{1,1} + \tfrac{1}{r'} & u_{2,1} + \tfrac{1}{r'} & \ldots & u_{m,1} + \tfrac{1}{r'} \\
			u_{1,2}     & u_{2,2}     & \ldots & u_{m,2} \\
			\vdots         & \vdots         & \ddots & \vdots     \\
			u_{1,m}     & u_{2,m}     & \ldots & u_{m,m}.
		\end{array}
		\right|}.
\end{align*}

Noticing that $-\tfrac{1}{r'}\zzz = (-\tfrac{1}{r'}, 0, \ldots, 0)^T$, we observe that each column of the matrix
is the projection onto $L=Q^\perp$ of the vector $(\uuu_i - (-\tfrac{1}{r'}\zzz))$, and the absolute value
of the determinant is thus $m! \Vol{m}{-\tfrac{1}{r'}\zzz, \uuu^{L}_1, \ldots, \uuu^{L}_m}$.
Combining all together, we get the final expression for the Jacobian:
\begin{align*}
	J &= (r^*)^{m-1} r^{m(n-1)} r' m! \Vol{m}{-\tfrac{1}{r'}\zzz, \uuu^{L}_1, \ldots, \uuu^{L}_m} \\
	  &= (r^*)^{m-2} r^{m(n-1) + 1} m! \Vol{m}{-\tfrac{r^*}{r}\zzz, \uuu^{L}_1, \ldots, \uuu^{L}_m}
		,
\end{align*}
which settles the $m = n-q$ case.

We now proceed to the $m < n-q$ case. Decompose each vector $\xxx_i$ as $\xxx_i = \xxx_i^Q + \xxx_i^{\perp}$,
where the first component lies in $Q$ and the second in $(n-q)$-dimensional subspace $Q^\perp$. We apply the
linear Blaschke--Petkantschin formula \eqref{eq:classic-bp} in $Q^\perp$ as follows:
\begin{align*}
	\MoveEqLeft \int\displaylimits_{\xxx \in (\Rspace^n)^{m}} f(\xxx) \diff \xxx = \!\!\!\!\!\!
	   \int\displaylimits_{\xxx_1^Q, \ldots, \xxx_m^Q \in Q}
		   \,\,
		   \int\displaylimits_{\xxx_1^{\perp}, \ldots, \xxx_m^{\perp} \in {Q^\perp}} \!\!\!\!\!\!
					f(\xxx_1^Q + \xxx_1^{\perp}, \ldots) \diff \xxx^{\perp} \diff \xxx^Q \\
		&= \!\!\!\!\!\!
		 \int\displaylimits_{\xxx_1^Q, \ldots, \xxx_m^Q \in Q}
			\int\displaylimits_{L \in \LGrass{m}{}(Q^\perp)}
			 \int\displaylimits_{\yyy_1, \ldots, \yyy_m \in L} \!\!\!\!\!\!
	    	f(\xxx_1^Q + \yyy_1, \ldots) \left[m! \Vol{m}{\yyy}\right]^{n-q - m} \diff \yyy \diff L \diff \xxx^Q.
\end{align*}
Now we apply Fubini's theorem and join integrations by $Q$ and $L$, and finish with applying the already proved result in $(q+m)$-dimensional space $L \oplus Q$.
Writing $\LGrass{m}{n-q}$ for the Grassmanian of $m$-subspaces inside $Q^\perp$, and $\yyy^{L}$ for the orthogonal projection of $\yyy$ onto $L$,
\begin{align*}
	\MoveEqLeft \int\displaylimits_{\xxx \in (\Rspace^n)^{m}} f(\xxx) \diff \xxx = 
		 \int\displaylimits_{L \in \LGrass{m}{n-q}}
			 \int\displaylimits_{\yyy_1, \ldots, \yyy_m \in L \oplus Q} \!\!\!\!\!\!
	    	f(\yyy_1, \ldots) \left[m! \Vol{m}{\yyy^{L}}\right]^{n-q - m} \diff \yyy \diff L \\
		&=
		 \int\displaylimits_{L \in \LGrass{m}{n-q}}
		   \int\displaylimits_{r \geq r_0}
        \int\displaylimits_{z \in S(L)} 
				 \int\displaylimits_{\uuu \in (S(L \oplus Q))^{m}}\!\!\!\!\!\!				
				  f(r^* z + r \uuu) \left[m! \Vol{m}{0, r^* z + r \uuu^L}\right]^{n-q - m} \times
				  \\
          &~~~~~~~~~~~~~~~~~\times r^{m(q+m-1) + 1} m! \Vol{m}{-\tfrac{\sqrt{r^2-r_0^2}}{r}\zzz, \uuu^{L}_1, \ldots, \uuu^{L}_m}
					\diff \uuu \diff z \diff r \diff L \\
		&= \int\displaylimits_{L \in \LGrass{m}{n-q}}
		   \int\displaylimits_{r \geq r_0}
        \int\displaylimits_{z \in S(L)} 
				 \int\displaylimits_{\uuu \in (S(L \oplus Q))^{m}} \!\!\!\!\!\!
				   f(r^* z + r \uuu)  \left[m! \Vol{m}{-\tfrac{r^*}{r}\zzz, \uuu^{L}_1, \ldots, \uuu^{L}_m}\right]^{n-q-m+1} \!\!\!\times \\
				   &~~~~~~~~~~~~~~~~~\times r^{m(q+m-1) + 1 + m(n-q-m)} (r^*)^{m-2}
					 \diff \uuu \diff z \diff r \diff L,
\end{align*}
the last equality being a consequence of the scaling of volume.
The bracket in the power of $r$ expands to $m(n-1)+1$, which completes the proof.

\bibliographystyle{abbrv}
\bibliography{bp}

\begin{thebibliography}{10}

\bibitem{Bla35}
W.~Blaschke.
\newblock {Integralgeometrie 2: Zu Ergebnissen von M.W. Crofton}.
\newblock {\em Bull. Math. Soc. Roum. Sci.}, 37:3--11, 1935.

\bibitem{EN17_1}
H.~Edelsbrunner and A.~Nikitenko.
\newblock Random inscribed polytopes have similar radius functions as
  poisson--delaunay mosaics.
\newblock {\em Annals of Applied Probability}, 28, 2017.

\bibitem{EN17_2}
H.~Edelsbrunner and A.~Nikitenko.
\newblock {Weighted Poisson--Delaunay mosaics}.
\newblock arXiv:1705.08735 [math.PR], 2017.

\bibitem{ENR16}
H.~Edelsbrunner, A.~Nikitenko, and M.~Reitzner.
\newblock {Expected sizes of Poisson--Delaunay mosaics and their discrete Morse
  functions}.
\newblock {\em Adv. Appl. Prob.}, 49:745--767, 2017.

\bibitem{Had57}
H.~Hadwiger.
\newblock {\em {Vorlesungen {\"u}ber Inhalt, Oberfl{\"a}che und
  Isoperimetrie}}, volume~93 of {\em Die Grundlehren der Mathematischen
  Wissenschaften}.
\newblock Springer-Verlag, Berlin, 1957.

\bibitem{Mil70}
R.~E. Miles.
\newblock {On the homogeneous planar Poisson point process}.
\newblock {\em Math. Biosci.}, 6:85--127, 1970.

\bibitem{Mil71}
R.~E. Miles.
\newblock {Isotropic random simplices}.
\newblock {\em Adv. Appl. Prob.}, 3:353--382, 1971.

\bibitem{Mol94}
J.~M{\o}ller.
\newblock {\em {Lectures on Random Voronoi Tessellations}}.
\newblock Springer New York, 1994.

\bibitem{Pet35}
B.~Petkantschin.
\newblock {Integralgeometrie 6. Zusammenh{\"a}nge zwischen den Dichten der
  linearen Unterr{\"a}ume imn- dimensionalen Raum}.
\newblock {\em Abhandlungen aus dem Mathematischen Seminar der Universit{\"a}t
  Hamburg}, 11(1):249--310, 1935.

\bibitem{ScWe08}
R.~Schneider and W.~Weil.
\newblock {\em {Stochastic and Integral Geometry}}.
\newblock Springer, Berlin, Germany, 2008.

\bibitem{Var36}
O.~Varga.
\newblock {Integralgeometrie 3. Croftons Formeln f{\"u}r den Raum}.
\newblock {\em Mathematische Zeitschrift}, 40(1):387--405, Dec 1936.

\end{thebibliography}

\end{document}